\documentclass[12pt]{article}
\def\date{December 2, 2015}
\usepackage{amsmath, amssymb, amsfonts, theorem, color, enumitem,graphicx}
\usepackage[margin=1in]{geometry}

\newtheorem{proposition}{Proposition}
\newtheorem{lemma}[proposition]{Lemma}

\newtheorem{theorem}[proposition]{Theorem}

\newcommand{\qed}{$\square$\bigskip}

{\theorembodyfont{\rmfamily}

}

\newcommand{\NN}{\mathbb{N}}

\begin{document}

\phantom{a}\vskip .25in
\centerline{{\large \bf  THE GY\H ORI-LOV\'ASZ THEOREM}%
\footnote{%
Partially supported by NSF under Grant No.~DMS-1202640.}}

\vskip.4in
\centerline{{\bf Alexander Hoyer}}

\centerline{and}

\medskip
\centerline{{\bf Robin Thomas}}
\medskip
\centerline{School of Mathematics}
\centerline{Georgia Institute of Technology}
\centerline{Atlanta, Georgia  30332-0160, USA}

%
%
%
\baselineskip 18pt

\vskip .75in

\noindent 
Our objective is to give a self-contained proof of the following beautiful theorem of Gy\H ori~\cite{gyo} and Lov\'asz~\cite{LovHomology},
conjectured and partially solved by Frank~\cite{FraPhD}.

\begin{theorem}
Let  $k\ge2$ be an integer, let $G$ be a $k$-connected graph on $n$ vertices, let $v_1,v_2,\ldots,v_k$ be distinct
vertices of~$G$, and let $n_1,n_2,\ldots,n_k$ be positive integers with $n_1+n_2+\cdots+n_k=n$.
Then $G$ has disjoint connected subgraphs $G_1,G_2,\ldots,G_k$ such that, for $i=1,2,\ldots,k$,
the graph $G_i$ has $n_i$ vertices and $v_i\in V(G_i)$.
\end{theorem}

\noindent 
The proof we give is  Gy\"ori's original proof, restated  using our terminology.
It clearly suffices to prove the following.

\begin{theorem}
\label{thm2}
Let  $k\ge2$ be an integer, let $G$ be a $k$-connected graph on $n$ vertices, let $v_1,v_2,\ldots,v_k$ be distinct
vertices of~$G$, and let $n_1,n_2,\ldots,n_k$ be positive integers with $n_1+n_2+\cdots+n_k<n$.
Let $G_1,G_2,\ldots,G_k$ be disjoint connected subgraphs of $G$ such  that, for $i=1,2,\ldots,k$,
the graph $G_i$ has $n_i$ vertices and $v_i\in V(G_i)$.
Then $G$ has disjoint connected subgraphs $G_1',G_2',\ldots,G_k'$ such that $v_i\in V(G_i')$ for $i=1,2,\ldots,k$,
the graph  $G_1'$ has  $n_1+1$ vertices and for $i=2,3,\ldots,k$
the graph $G_i'$ has $n_i$ vertices.
\end{theorem}

For the proof of Theorem~\ref{thm2} we will use terminology inspired by hydrology (the second author's father would have been pleased).
Certain vertices will act as ``dams" by blocking other vertices from the rest of a subgraph of $G$, thus creating a ``reservoir".
A sequence of dams will be called a ``cascade".

To define these notions precisely let $G_1,G_2,\ldots,G_k$ be as in Theorem~\ref{thm2} and 
let $i=2,3,\ldots,k$.
For a vertex $v\in V(G_i)$ we define the \textbf{reservoir} of $v$, denoted by $R(v)$, to be the set of all vertices in~$G_i$ which 
are connected to $v_i$ by a path in $G_i\backslash v$. Note that $v\notin R(v)$ and also $R(v_i)=\emptyset$.
By a \textbf{cascade} in $G_i$ we mean a (possibly null) sequence $w_1, w_2,\ldots, w_m$ of distinct vertices in $G_i\backslash v_i$
such that $w_{j+1}\notin R(w_j)$ for $j=1,\ldots,m-1$.
Thus $w_j$ separates $w_{j-1}$ from $w_{j+1}$ in~$G_i$ for every $j=1,\ldots,m-1$, where $w_0$ means $v_i$.
By a \textbf{configuration} we mean a choice of subgraphs $G_1,G_2,\ldots,G_k$ as in Theorem~\ref{thm2}
and exactly one cascade in each $G_i$ for $i=2,3,\ldots,k$.
By a \textbf{cascade vertex} we mean a vertex belonging to one of the cascades in the configuration.
We define  the \textbf{rank} of some cascade vertices recursively as follows. 
Let $w\in V(G_i)$ be a cascade vertex. If $w$ has a neighbor in $G_1$, then we define the rank of $w$  to be $1$. 
Otherwise, its rank is the least integer $k\ge2$ such that there is a cascade vertex $w'\in V(G_j)$, for some $j\in\{2,3,\ldots,k\}-\{i\}$, so that $w$ has a neighbor in $R(w')$ and $w'$ has rank $k-1$. If there is no such neighbor, then the rank of $w$ is undefined.
For an integer $r\ge1$, let $\rho_r$ denote the total number of vertices belonging to $R(w)$ for some cascade vertex $w$ of rank $r$.
A configuration is \textbf{valid} if each cascade vertex has well-defined rank and this rank is strictly increasing within a cascade. That is, for each cascade $w_1,w_2,\ldots,w_m$ and integers $1\le i<j\le m$ the rank of $w_i$ is strictly smaller than the rank of $w_j$.
 Note that a valid configuration exists trivially by taking each cascade to be the null sequence.
For an integer $r\ge1$  a valid configuration is \textbf{$r$-optimal} if, among all valid configurations, 
it maximizes $\rho_1$, 
subject to that it maximizes $\rho_2$, and so on, up to maximizing $\rho_r$. If a valid configuration is $r$-optimal for all $r\ge1$, we simply say it is \textbf{optimal}.

Finally, we define $S:=V(G)-V(G_1)-V(G_2)-\cdots-V(G_k)$. This is nonempty in the setup of Theorem 2. We say that
a \textbf{bridge} is an edge with one end in $S$ and the other end in the reservoir 
of a cascade vertex. In a valid configuration, the \textbf{rank} of the bridge is the minimum rank of all cascade vertices $w$ where the bridge has an end in $R(w)$.

\begin{figure}[h]
\centering
\includegraphics[scale=0.8]{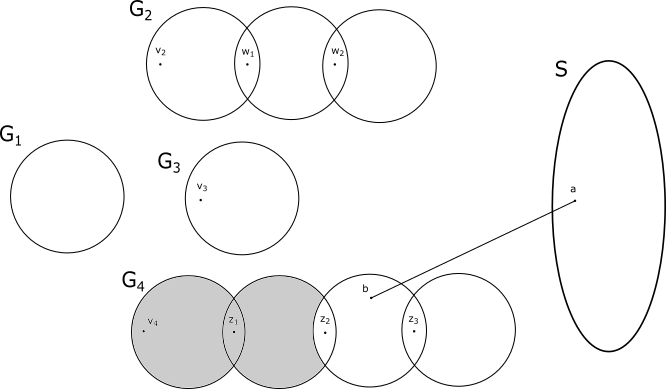}
\caption{An example of a configuration. $w_1, w_2, z_1, z_2,$ and $z_3$ are cascade vertices. $R(z_2)$ is shaded. The edge $ab$ is a bridge, and its rank is the rank of $z_3$.}
\label{fig1}
\end{figure}

These concepts are illustrated in Figure~\ref{fig1}.

\begin{lemma}
\label{lma1}
If there is an optimal configuration containing a bridge, then the conclusion of Theorem~\ref{thm2} holds.
\end{lemma}

\noindent 
{\bf Proof.}
Suppose there is an optimal configuration containing a bridge. Then for some $r\in\NN$ we can find a configuration which is $r$-optimal containing a bridge of rank $r$. Choose the configuration and bridge so that $r$ is minimal. Denote the endpoints of the bridge as $a\in S$ and 
$b\in R(w)\subseteq V(G_i)$, where $w$ is a cascade vertex of rank $r$.

Suppose $w$ separates $G_i$. Since we have a valid configuration, any cascade vertices in 
$V(G_i)-R(w)-\{w\}$
must have rank greater than $r$. Choose any nonseparating vertex from this set, say $u$. We make a new valid configuration in the following way. Move $u$ to $S$ and $a$ to $G_i$. Leave the cascades the same with one exception:
remove all cascade vertices in $V(G_i)-R(w)-\{w\}$ and all cascade vertices whose rank becomes undefined.
Note that any cascade vertices affected by this action have rank greater than $r$. Now our new configuration is valid, increased the size of $R(w)$, and did not change any other reservoirs of rank at most $r$. This contradicts $r$-optimality.

So, continue under the assumption that $w$ does not separate $G_i$. 
If $r=1$, choose $G_1':=G_1+w$, the graph obtained from $G_1$ by adding the vertex $w$ and all edges from $w$ 
to $G_1$, $G_i':=(G_i+a)\backslash w$, and leave all other $G_j$'s unchanged.
Then these graphs  satisfy the conclusion of Theorem~\ref{thm2}, as desired.

If $r>1$, then $w$ has a neighbor in some $R(w')$ with $\text{rank}(w')=r-1$. As before, we make a new valid configuration by moving $w$ to $S$ and $a$ to $G_i$. Keep the cascades the same as before, except terminate $w$'s former cascade just before $w$ and exclude any cascade vertices whose rank has become undefined. Though we may have lost several reservoirs of rank $r$ and above, the new configuration is still $(r-1)$-optimal. Also, the edge connecting $w$ to its neighbor in $R(w')$ is now a rank $r-1$ bridge. This contradicts the minimality of $r$, so the proof of Lemma~\ref{lma1} is complete.~\qed

\begin{lemma}
\label{lma3}
Suppose there is an optimal configuration with an edge $ab$ such that:
\begin{enumerate}
\item Either $a\in V(G_1)$ or $a$ is in a reservoir, and
\item $b\in V(G_i)$ for some $i\in\{2,3,\ldots,k\}$, $b\neq v_i$, and $b$ is not in a reservoir.
\end{enumerate}
Then the cascade of $G_i$ is not null and $b$ is the last vertex in the cascade.
\end{lemma}

\noindent 
{\bf Proof.} Suppose there is such an edge in an optimal configuration and $b$ is not the last vertex in the cascade of $G_i$. Denote the cascade of $G_i$ by $w_1,\ldots, w_m$ (which a priori could be null). Since $b$ is not in a reservoir and is not the last cascade vertex, we know that $b$ is not a cascade vertex. Then make a new configuration by including $b$ at the end of $G_i$'s cascade. By condition 1, $b$ has well-defined rank. If this rank is larger than all other ranks in the cascade (including the case where the former cascade is null), then we have a valid configuration and have contradicted optimality by adding a new reservoir (which is nonempty since $v_i\in R(b)$) without changing anything else.

So, the former cascade is not null.
Let  $\text{rank}(b)=r$ and let $j\ge0$ be the integer such that $j=0$ if $r\leq\text{rank}(w_1)$ and 
 $\text{rank}(w_j)<r\leq\text{rank}(w_{j+1})$ otherwise.
We make a second adjustment by excluding the vertices $w_{j+1},w_{j+2},\ldots,w_m$ from the cascade and adding $b$
to it. Now the configuration is clearly valid, but it is unclear whether optimality has been contradicted. But notice that every vertex which used 
to belong to $R(w_{j+1})\cup R(w_{j+2})\cup\cdots\cup R(w_m)$ now belongs to $R(b)$, and also $R(b)$ contains $w_m$ which was not in any reservoir previously. Thus, we have strictly increased the size of rank $r$ reservoirs without affecting any lower rank reservoirs. This contradicts optimality, so the proof of Lemma~\ref{lma3} is complete.~\qed

\noindent 
{\bf Proof of Theorem~\ref{thm2}.}
Using our lemmas, we can assume we have an optimal configuration which does not contain any bridges and where any edges as in Lemma~\ref{lma3} are at the end of their cascades. Consider the set containing the last vertex in each non-null cascade and the $v_i$ corresponding to each null cascade. This is a cut of size $k-1$, separating $G_1$ and the reservoirs from the rest of the graph, including $S$. 
This contradicts $k$-connectivity, and the proof is complete.~\qed

\end{document}